\documentclass[twoside,reqno]{amsart}
\usepackage[utf8]{inputenc}
\usepackage{bm}
\usepackage{subcaption}

\usepackage{color}
\usepackage{amssymb}
\usepackage{graphicx}
\usepackage{float}

\newtheorem{theorem}{Theorem}[section]
\newtheorem{lemma}[theorem]{Lemma}

\numberwithin{equation}{section}

\DeclareMathOperator{\diag}{diag}

\renewcommand{\kappa}{\varkappa}
\newcommand{\be}{\begin{equation}}
\newcommand{\ee}{\end{equation}}

\newcommand{\bq}{\begin{eqnarray}}
\newcommand{\eq}{\end{eqnarray}}

\newcommand{\ba}{\begin{array}}
\newcommand{\ea}{\end{array}}

\newcommand{\al}{\alpha}

\newcommand{\twotwo}[4]{\left(\begin{array}{cc}#1&#2\\&\\#3&#4\end{array}\right)}

\title{Analytic solutions to nonlinear ODEs via spectral power series}
\author{Estelle Basor and Rebecca Morrison}
\date{\today}

\begin{document}

\begin{abstract}
    Solutions to most nonlinear ordinary differential equations (ODEs) rely on numerical solvers, but this gives little insight into the nature of the trajectories and is relatively expensive to compute. In this paper, we derive analytic solutions to a class of nonlinear, homogeneous ODEs with linear and quadratic terms on the right-hand side. We formulate a power series expansion of each state variable, whose function depends on the eigenvalues of the linearized system, and solve for the coefficients using some linear algebra and combinatorics. Various experiments exhibit quickly decaying coefficients, such that a good approximation to the true solution consists of just a few terms.\\
    
\end{abstract}

\maketitle
\section{Introduction}
Only a few systems of coupled ordinary differential equations (ODEs) admit known closed-form solutions, in particular, linear systems and the 1-dimensional logistic equation. When a closed-form or analytic solution  is not available, mathematicians, scientists, students, and engineers are left instead with numerical solutions: These are relatively slow to compute (usually requiring integration), difficult to analyze (how quickly the solution approaches equilibrium, how good are approximations to the true solution, or how the system might change with different inputs), require the introduction of new parameters (grid size or polynomial degree, for example), and require thorough verification tests.

In this work, we provide closed-form solutions for a class of first-order nonlinear ODEs with linear and quadratic terms on the right-hand side (RHS). Such systems are mainstays in diverse fields including the common SEI or SEIR-type models of epidemiology~\cite{kuhl2021computational}; chemical reaction mechanisms in combustion and other chemical kinetics~\cite{williams2018combustion, espenson1995chemical}; and predator-prey or Lotka-Volterra models in mathematical biology~\cite{strogatz2018nonlinear}. In all of these examples, the rates of change of different state variables---epidemiological/chemical/biological species---depend on some linear growth rates and pairwise (quadratic) interaction rates.

Our solutions describe the behavior of the system for all $t >t_0$ where $t_0$ is some fixed value, known exactly in the one-dimensional case, and with a bound provided in the higher-dimensional cases. We assume the eigenvalues of the linearized system at the equilibrium are real and negative and independent over the field of rationals. The method is based on power series expansions along with some linear algebra and a bit of combinatorics; we refer to it as the Spectral Power Series (SPS) method. Compared to existing work, SPS overlaps with approaches using the eigendecomposition of the Koopman operator~\cite{mauroy2020koopman}. On one hand, SPS is so far limited to linear and quadratic RHS, while Koopman operator theory applies more generally. On the other hand, in contrast with Koopman methods, here we generate explicit functions and coefficients comprising the infinite series. We also note the work of~\cite{kaniadakis2022novel}, in which the two-dimensional Lotka-Volterra equations are reframed with a new Hamiltonian formalism that does admit an analytic solution but replaces the competition terms with other functions from a specified class.

Our motivation to find such solutions came from the study of reduced models. In applications, certain perturbation methods are useful for finding subsets of solutions for complicated models with incomplete data. Our goal was to understand why sometimes such techniques work and why sometimes they fail, but reliance on numerical solutions stalled much progress. More will be said about this topic at the end of the paper. 

As a proof of concept, we demonstrate the method for the logistic equation: 
\begin{equation}\dot x = rx\left(1-\frac{x}{k}\right),\label{eq:logistic} \end{equation}
which has the known solution
\[ x(t) = \frac{k}{1 + Ae^{-rt}}\] where $A = (k-x_0)/x_0$, initial condition $x_0 = x(0)$, and equilibrium $x_{eq} = k$.
We will recover this solution with a power series expansion of $x$. That is, first assume that $x(t) = \sum_n \al^{*}_n (g(t))^n$ for some unknown function $g$. We substitute this into the differential equation:
\[\sum_n \alpha^*_n n g^{n-1}(t) g' (t) = r\sum_n \alpha^*_ng^n(t) - \frac{r}{k}\left(\sum_n \alpha^*_n g^n(t) \right)\left( \sum_m \alpha^*_m g^m(t) \right).\]
We first notice that matching coefficients for $n=0$ recovers $\alpha_0 = k$, i.e., the constant term is equal to the equilibrium.  For $n=1$, we find $g'/g = -r$, so $\ln g = -rt + \text{const.}$, and thus $g(t) = ce^{-rt}$. Let $\alpha_n = \alpha_n^* c^{n}$, so that
\[x(t) = \sum_n \alpha_n e^{-rnt}.\]

Higher-order terms yield $\alpha_n = \alpha_1^n/k^{n-1}, n \geq 2$.
That is, $\al_1$ is so far left undetermined, and all higher-order coefficients are written in monomials in terms of $\al_1$. So the solution becomes
\[x(t) =  k\sum_{n\geq 0} \left(\frac{\al_1}{k}\right)^n e^{-r n t}.\] Given some fixed $\al_{1}$ this series converges for large enough $t$ and is
\[ x(t) = \frac{k^{2}}{ k - \al_{1}e^{-rt}}. \] 
Suppose that $|\al_1/k| < 1$. Then
this solution is valid for all $t \geq 0,$ and $\al_1$ is given by the initial condition:
\[x_0 = \frac{k^2}{k-\al_1}\]
which implies $\al_1 = k(\frac{x_0 - k}{x_0})$. 
Finally, we plug this into the solution to find
\begin{align*}
x(t) = \frac{k}{1 + Ae^{-rt}},
\end{align*}
where A is the same as defined above.

If we parametrize $\al_1$ as $k(\frac{x_0 - k}{x_0})$, note the above condition is satisfied when $x_{0} > \frac{1}{2}k.$ However, even if $0 < x_0 \leq \frac{1}{2}k$, the series solution extends analytically for all $t > 0.$ In general the series converges for 
\[t> t_0 = \frac{1}{r}(\log| 1- k/x_0 |).\]
It is also worth pointing out that \[ \al_1 = \lim_{t \rightarrow \infty} (x(t) - k)e^{rt},\] that is, the solution at infinity retains knowledge about initial conditions. 

The logistic solution above can of course be found using simply a separation of variables technique. But the point is it can also be easily found using a series.  For the logistic equation, the functions $e^{-rnt}$ and coefficients $\al_n = \left(\frac{x_{0} - k}{x_{0}}\right)^{n}$ above correspond exactly to a set of eigenvalues and eigenfunctions, respectively, of the Koopman operator flow~\cite{li2023adaptive}. In particular, the choice to parametrize $\al_{1}/k$ as $\frac{x_0 - k}{x_0}$ yields a Koopman expansion. 

This same SPS technique works for systems that are analogous to the logistic equation. We show how this is possible in the next sections and investigate the domains of convergence. For these systems as well, the terms of the series overlap with the Koopman expansion found in~\cite{mezic2019}. Unlike most of the Koopman literature, we show that the series converges and find estimates for the coefficients and for the domains of convergence. Our methods are linear algebra-based and use a constructive recursive algorithm to find the coefficients. We will say more later in the paper about connections to the Koopman expansions. 

\section{$2\times2$ systems}\label{sec:2by2}
In this section, we address the $2 \times 2$ case in detail. It is a bit easier to track the steps for this case and then it turns out that the higher dimensional case only needs minor modifications in the argument. 

We consider 2-dimensional systems of the form
\begin{align*} \dot{x_1} &= b_1x_1 + a_{11} x_1^{2} + a_{12} x_1 x_2 \\
\dot{x_2} &=  b_2 x_2  + a_{21} x_1x_2 + a_{22} x_2^2.\end{align*}
Let $$ \bm A = \twotwo{a_{11}}{a_{12}}{a_{21}}{a_{22}}.$$ The equilibrium of the system is found by computing
$\bm c = \bm A^{-1}(-\bm b)$ where $$\bm b = \left(\begin{array}{c}b_1\\b_2\end{array}\right)
\mbox{and }
\bm c = \left(\begin{array}{c}c_1\\c_2\end{array}\right).$$

The eigenvalues of the linearized system near equilibrium are given by the eigenvalues of the matrix
$$ \text{diag}(\bm c) \bm A =  \twotwo{c_1a_{11}}{c_1a_{12}}{c_2a_{21}}{c_2a_{22}}.$$
 To see this, let $\bm x = \bm c+ \bm y$ around equilibrium, so that $\dot {\bm x} = \dot{\bm y}$ and
\begin{align*}
    \dot{\bm y} &= \diag(\bm y+ \bm c) \left( \bm b + \bm A (\bm c+ \bm y)\right)\\
    &= \diag(\bm c+ \bm y) \left( \bm A \bm y\right)
    \end{align*}
    and the linear part is $\dot{\bm y} \approx \diag(\bm c) \bm A \bm y$.
The eigenvalues, which we will denote by $\lambda_1$ and $\lambda_2$ are given by 
$$ \frac{(c_1a_{11} +c_2a_{22}) \pm \sqrt{(c_1a_{11} - c_2a_{22})^2 +4c_1c_2a_{12}a_{21}}}{2}.$$

We now assume that both eigenvalues are real, negative, and that they are independent over the field of rationals, and that our solutions near equilibrium are of the form
\begin{align*} x_1(t) &=  \sum_{n_{1},n_{2} \geq 0} \al_{1}^{n_{1},n_{2}} e^{(n_{1}\lambda_1 + n_{2}\lambda_2)t} \\
 x_2(t) &= \sum_{n_{1},n_{2} \geq 0} \al_{2}^{n_{1},n_{2}} e^{(n_{1}\lambda_1 + n_{2}\lambda_2)t}.\end{align*}
(From here on, we use subscripts for variable or eigenvalue indices, and superscripts for the expansion indices.)
Substituting these solutions into the differential equations and equating like terms, we find that as expected 
$\al_{1}^{0,0} = c_1, \,\,\,\,\al_{2}^{00} = c_2.$ 
For the cases when $n_{1} = 1, n_{2} = 0,$ and $n_{1} = 0, n_{2} = 1,$ the equations leave $\al_1^{1,0}$ and $\al_1^{0,1}$ undetermined, $\al_{2}^{1,0}$ and $\al_{2}^{0,1}$ multiples of the undetermined coefficients, but yield values for $\lambda_1$ and $\lambda_2$ that agree with the ones found above.
(To see this, match coefficients for $e^{\lambda_1 t}$ or $e^{\lambda_2 t}$ terms. Both result in a quadratic equation for $\lambda_{1,2}$.)

Returning to the coefficients, 
\[\al_{2}^{1,0} = \al_{1}^{1,0}a_{21}c_2/( \lambda_1 - a_{22}c_2)\] and
\[\al_{2}^{0,1} = \al_{1}^{0,1}a_{12}c_2/( \lambda_2 - a_{11}c_1).\] 

For all other values we see that $\al_{1}^{n_{1},n_{2}}$ and $\al_{2}^{n_{1},n_{2}}$ must satisfy
\begin{align} (n_{1} \lambda_1\ + n_{2}\lambda_2 - a_{11} c_1)\al_{1}^{n_{1},n_{2}} - a_{12} c_1 \al_{2}^{n_{1},n_{2}} &= a_{11} S_{1,1}^{n_{1},n_{2}}+ a_{12} S_{1,2}^{n_{1},n_{2}} \nonumber \\
 -a_{21}c_2\al_{1}^{n_{1},n_{2}} + (n_{1} \lambda_1\ + n_{2}\lambda_2 - a_{22} c_2)\al_{2}^{n_{1},n_{2}} &= a_{21} S_{1,2}^{n_{1},n_{2}} + a_{22} S_{2,2}^{n_{1},n_{2}} \label{eq:2sys}\end{align}
where
\[ S_{1,1}^{n_{1},n_{2}} = \sum_{(i,j) \neq (0,0),(n_{1},n_{2})} \al_{1}^{i,j}\al_{1}^{n_{1}-i,n_{2}-j} \]
\[S_{1,2}^{n_{1},n_{2}} = \sum_{(i,j) \neq (0,0),(n_{1},n_{2})} \al_{1}^{i,j}\al_{2}^{n_{1}-i,n_{2}-j}\]
and
\[ S_{2,2}^{n_{1},n_{2}} = \sum_{(i,j) \neq (0,0),(n_{1},n_{2})} \al_{2}^{i,j}\al_{2}^{n_{1}-i,n_{2}-j} .\]

Now let $\bm n = \langle n_{1},n_{2} \rangle, \bm \lambda = \langle \lambda_{1},\lambda _{2} \rangle,$  $\bm \al^{\bm n} = \langle \al_{1}^{n_{1},n_{2}}, \al_{2}^{n_{1},n_{2}} \rangle,$ 
\[ \bm S^{\bm n} = \twotwo{S_{1,1}^{n_{1},n_{2}}}{S_{1,2}^{n_{1},n_{2}}}{S_{1,2}^{n_{1},n_{2}}}{S_{2,2}^{n_{1},n_{2}}},\] 
and 
\[ \bm s^{\bm n} = \langle \bm e_{1}\cdot \bm A \bm S_{\bm n}\bm e_{1}, \bm e_{2}\cdot \bm A \bm S_{\bm n}\bm e_{2} \rangle \]with
$\bm e_{1}, \bm e_{2}$ the standard unit basis vectors. In other words, the vector $\bm s^{\bm n}$ is the diagonal of $\bm A \bm S^{\bm n}.$

Our system (\ref{eq:2sys}) is equivalent to 
\[(\bm n \cdot \bm \lambda) \left ( \bm I - \frac{\diag({\bm c}) \bm A}{\bm n \cdot \bm \lambda } \right) \bm{\al}^{\bm n} = \bm s^{\bm n}\]

or 

\[ \bm \al^{\bm n} = \frac{1}{\bm n \cdot \bm \lambda} \left( \bm I - \frac{\diag({\bm c}) \bm A}{\bm n \cdot \bm \lambda }\right)^{-1}\bm s^{\bm n}.\]

Note that the matrix $\left( I - \frac{\diag({\bm c}) A}{\bm n \cdot \bm \lambda }\right)$ is always invertible if $\bm n$ is not $(1,0)$ or $(0,1).$
%
%
%

Our goal now is to estimate the terms in the last equation so that we can get a good estimate of $\bm \al^{\bm n}.$

We already assumed that both $\lambda_1$ and $\lambda_2$ are real,  negative, and independent over the rationals. We now also assume that $|\lambda_1| < |\lambda_2|$ without any loss of generality.

\begin{lemma}
Let $N = n_{1}+ n_{2}.$ Then
\[ |\bm n \cdot \bm \lambda| \geq |\lambda_{1}|N.\]
\end{lemma}
\proof
This follows easily since all the eigenvalues are negative.
\endproof
\begin{lemma}\label{S.est}Let $N = n_{1}+ n_{2}.$
Suppose not both $i$ and $j$ are zero and suppose for all other $i$ and $j$ with $i + j <  N$ there exists a constant K such that 
\[|\al_{1}^{i,j}| \leq \frac{K^{i+j}}{(i+1)(j+1)},\,\,\,\, |\al_{2}^{i,j}| \leq \frac{K^{i+j}}{(i+1)(j+1)}.\]
Then 
    \[|S_{1,1}^{n_{1},n_{2}}| \leq \frac{4 K^{n_{1}+n_{2}}((\log (n_{1}+1))+1)((\log(n_{2}+1)+1)}{(n_{1}+1)(n_{2}+1)} .\]

\end{lemma}
\proof
We have that
\[ S_{1,1}^{n_{1},n_{2}} = \sum_{(i,j) \neq (0,0),(n_{1},n_{2})} \al_{1}^{i,j}\al_{1}^{n_{1}-i,n_{2}-j} \]
is bounded by 
\[\sum_{(i,j) \neq (0,0),(n_{1},n_{2})} \frac{K^{i+j}}{(i+1)(j+1)}\frac{K^{n_{1}-i+ n_{2}-j}}{(n_{1}-i+1)(n_{2}-j+1)}\]
or
\[K^{n_{1}+n_{2}}\sum_{i=0}^{n_{1}} \frac{1}{(i+1)(n_{1}-i+1)} \sum_{j=0}^{n_{2}}\frac{1}{(j+1)(n_{2}-j+1)}\]
Separating the denominators into partial fractions and using an estimate for $\log n$ then yields the result.
\endproof
\begin{lemma}
 Assuming the conditions in the previous lemma, the same estimates hold for  
 and $S_{1,2}^{n_{1},n_{2}}$ and $S_{2,2}^{n_{1},n_{2}}.$
\end{lemma}
\begin{lemma}
Suppose the coefficients of the two power series 
\[ \sum_{n_{1},n_{2} \geq 0} \al_{1}^{n_{1},n_{2}} e^{(n_{1}\lambda_1 + n_{2}\lambda_2)t} \]
\[ \sum_{n_{1},n_{2} \geq 0} \al_{2}^{n_{1},n_{2}} e^{(n_{1}\lambda_1 + n_{2}\lambda_2)t} \]
satisfy the conditions
\[ |\al_{1}^{n_{1},n_{2}}| < \frac{C^{n_{1}+n_{2}}}{(n_{1}+1)(n_{2}+1)}, \,\,\, \,\,|\al_{2}^{n_{1},n_{2}}| < \frac{C^{n_{1}+n_{2}}}{(n_{1}+1)(n_{2}+1)} \]
for some positive constant $C$. Then there exists a $t_0$ such that for all $t > t_0$ the series converge. 
\end{lemma}

Note that once the above is verified, then term by term differentiation of the series yield series that converge uniformly on the same interval and thus the series are  the actual solutions.

\proof
The proof is the same for both series, so we consider
\[\sum_{n_{1},n_{2} \geq 0} \al_{1}^{n_{1},n_{2}} e^{(n_{1}\lambda_1 + n_{2}\lambda_2)t}  .\]
This sum is at most 
\[\sum_{n_{1},n_{2} \geq 0} \frac{C^{n_{1}+n_{2}}}{(n_{1}+1)(n_{2}+1)} e^{(n_{1} \lambda_1 + n_{2} \lambda_2)t}\]
which is bounded by
\[\sum_{n_{1} \geq 0} C^{n_{1}} e^{(n_{1}\lambda_1)t} \,\,\sum_{n_{2} \geq 0} C^{n_{2}}e^{(n_{2}\lambda_1)t}\]
Thus if $t_0$ is chosen so that  $|C e^{ \lambda_1 t_0}| <1$ and $|C e^{ \lambda_2 t_0}| <1$
we have convergence.
\endproof

%

Our last step is to show how one can guarantee the estimates in the last lemma. In the lemma we 
use a few facts about matrix norms. 

The standard operator norm for a matrix $A = (a_{i,j}),$  is defined to be \\
$\|A\|=\sup_{||x|| =1}\|A(x)\|,$
while the Hilbert-Schmidt norm is defined by \\$\|A\|_2 = (\sum_{i,j}|a_{ij}|^2)^{1/2}.$
It is well known that 
$\|A\| \leq \|A\|_2$ and that if $A = BC,$ then $\|A\|_2 \leq \|B\|\|C\|_2.$ (The same inequalities are true for bounded operators.) Details can be found in \cite{GK}.

\begin{theorem}
There exists a constant $K$ such that the coefficients $\al_{1}^{n_{1},n_{2}}$ and $\al_{2}^{n_{1},n_{2}}$ satisfy
\[ |\al_{1}^{n_{1},n_{2}}| < \frac{K^{n_{1}+n_{2}}}{(n_{1}+1)(n_{2}+1)}, \,\,\, \,\,|\al_{2}^{n_{1},n_{2}}| < \frac{K^{n_{1}+n_{2}}}{(n_{1}+1)(n_{2}+1)}.\]
\end{theorem}
\proof

Recall that 

\[ \bm \al^{\bm n} = \frac{1}{\bm n \cdot \bm \lambda} \left( \bm I - \frac{\diag{(\bm c)} \bm A}{\bm n \cdot \bm \lambda }\right)^{-1} \bm s^{\bm n}\]
We first note that a matrix of the form $\bm I + \bm B$ where $\bm B$ has operator norm $||\bm B|| \leq \delta < 1$ has an inverse of the form 
$$\bm I + \bm C$$ where the norm of $\bm C$ is at most $\frac{\delta}{1 - \delta}.$

Now choose $N_{0}$ such that for $N \geq N_{0}$ the operator norm 
\[ \Big\| \frac{\diag{(\bm c)} \bm A}{ \bm n \cdot \bm \lambda }\Big\| \,\,\,\leq \,\,\,\Big\|\frac{\diag{(\bm c)} \bm A}{N\lambda_{1} } \Big\|\] is at most $1/2.$ Thus the operator norm of 
\[\left( \bm I - \frac{\diag{(\bm c)} \bm A}{\bm n \cdot \bm \lambda }\right)^{-1}\]
will be at most $ 2.$

Also, there exists an $N_{1}$ such that for $N > N_{1},$ 
\[   \frac{8( \log (N+1)+1)^{2} \|\bm A\|}{N|\lambda_1 |} < 1/2 .\]

Now let $N_{2} = \max\{N_{0},N_{1}\}.$


 For pair $i,j$ such that $i + j \leq N_2$ we compute the quantity $((i+1)(j+1)|a_{i,j}|)^{\frac{1}{i+j}}.$ There are a finite number of these and we define 
$K_1$ to be the maximum value of these. Thus for a finite number of such coefficients, 
we have that 
\[ |\al_{1}^{i,j} | < \frac{K_1^{i+j}}{(i+1)(j+1)}.\]

We repeat the process for the $\al_{2}^{i,j}$ and then choose the largest: $K = \max\{K_1,K_2\}$. We are finally in a position to use an induction argument. 

We let $n_{1}+ n_{2} = N_2 + 1.$ For all $i,j$ with $i+j \leq N_2$ we have the satisfied the condition
\[ |\al_1^{i,j} | < \frac{K^{i+j}}{(i+1)(j+1)},\,\,\,\,  |\al_2^{i,j} | < \frac{K^{i+j}}{(i+1)(j+1)}.\]

Thus by Lemma~\ref{S.est}, for any $n_{1}+ n_{2}$ with $n_{1}+ n_{2}= N_2 + 1,$ we have that 
$S_{i,j}^{n_{1}, n_{2}}$  ($i,j = 1,2$) is bounded by
\[\frac{4 K^{n_{1}+ n_{2}}(\log (n_{1}+1)+1)(\log(n_{2}+1)+1)}{(n_{1}+1)(n_{2}+1)}.\] 

This means the operator norm of the matrix $\bm S^{\bm n}$ is bounded by 
\[ \frac{8 K^{n_{1}+ n_{2}}((\log (N_{2}+2))+1)^{2}}{(n_{1}+1)(n_{2}+1)}.\] 
Note that the Hilbert-Schmidt norm of $\bm S^{\bm n}, \|\bm S^{\bm n}\|_{2}$, is bounded by the same value 
and thus 
\[ \| \bm s^{\bm n} \| \leq \|\bm A \bm S^{\bm n}\|_{2} \leq \|\bm A\| \|\bm S^{\bm n}\|_{2}.\]
Putting everything together we have for $\bm{n}  = \langle n_{1}, n_{2} \rangle,$ and $ n_{1}+ n_{2 } = N_{2} + 1$,
\[ \bm \al^{\bm n} = \frac{1}{\bm n \cdot \bm \lambda} \left( I - \frac{\diag({\bm c}) \bm A}{\bm n \cdot \bm \lambda }\right)^{-1}\bm s^{\bm n}\]
so that
\begin{align*}
\| \bm \al^{\bm n} \| &\leq \frac{1}{(N_{2} + 1)|\lambda_{1}|)} \Big \|\left( I - \frac{\diag({\bm c}) \bm A}{\bm n \cdot \bm \lambda }\right)^{-1}\Big \| \|\bm A\| \|\bm S^{\bm n}\|_{2} \\
&\leq \frac{2 \|\bm A\|}{(N_{2} + 1)|\lambda_{1}|)}   \frac{8 K^{n_{1}+ n_{2}}((\log (N_{2}+2))+1)^{2}}{(n_{1}+1)(n_{2}+1)} \\
& \leq \frac{ K^{n_{1}+ n_{2}}}{(n_{1}+1)(n_{2}+1)}\frac{2 \|\bm A\| 8(\log (N_{2}+2))+1)^{2}}{(N_{2} + 1)|\lambda_{1}|)}  
\end{align*} 
which implies
\[ \| \bm \al^{\bm n} \| \leq \frac{ K^{n_{1}+ n_{2}}}{(n_{1}+1)(n_{2}+1)}\]
or the desired result.
\endproof

One should note that the conditions that we imposed on the coefficients of our series were designed to make the analysis possible. There may be other approaches with less restrictive bounds that show the series converge. 

The other thing to note is that the choice of $1/2$ to find the $N_0$ and $N_1$ in the above argument can be replaced by $\delta < 1$ for $N_0$ and $1 - \delta$ for $N_1$ if it makes the choice of $N_2$ smaller. 

As we will see with some of the examples, the coefficients with choices of reasonable constants seem to decay quickly. This is expected because of the following result.

\begin{theorem}
Given some choice of $\al_1^{1,0}$ and $\al_1^{0,1}$ all other $\al_i^{n_1,n_2}, i = 1,2$ have the form 
\[ c_{\bm n} (\al_1^{1,0})^{n_{1}} (\al_1^{0,1})^{n_{2}} \] where the constant 
$c_{\bm n}$ only depends on all the other data in the system. 
\end{theorem}
\proof
This is straightforward to see from the formula for the coefficients
\[ \bm \al^{\bm n} = \frac{1}{\bm n \cdot \bm \lambda} \left( \bm I - \frac{\diag{(\bm c)} \bm A}{\bm n \cdot \bm \lambda }\right)^{-1}\bm s^{\bm n}\]
and the form of the terms for the $ S_{i,j}^{\bm n}.$
\endproof

We end this section with an example that estimates the value of $K$ and hence the value of $t_0$ in our main theorem. Other examples will be given in Section 4. 

For this example we let 
\[ \bm A = \twotwo{-1}{-0.3}{-0.1}{-1},\,\,\,\bm b = \left(\begin{array}{c}2.45\\1.7\end{array}\right),\,\,\, \bm c = \left(\begin{array}{c}2\\1.5\end{array}\right),\]with initial condition
$(3,3).$
The eigenvalues of $\diag (\bm c )\bm A$ are $\,\sim~-1.359$ and $\,\sim~-2.140 $. Using an estimate for the operator norm of $\bm A$ and some help from Mathematica we found that $N_2$ is of the order of~$ \sim 320.$ This yielded a $K$ value of less than $5.83$ and we find that we have convergence for $t > t_0 = 1.3.$

In addition, for this example, our eigenvalues are not too far apart in value, so we can find appropriate 
$\al_{1}^{1,0}$ and $\al_{1}^{0,1}$ by using limits at infinity. In fact
\[\al_{1}^{1,0} = \lim_{t \rightarrow \infty} (x_{1}(t) - 2)e^{-\lambda_1 t}\]
and
\[\al_{1}^{0,1} = \lim_{t \rightarrow \infty} (x_{1}(t) - 2 - \al_{1}^{1,0}e^{\lambda_1 t})e^{-\lambda_2 t}.\]
For our example, 
$\al_{1}^{1,0} \sim -0.45$ and
$\al_{1}^{0,1} \sim 0.91.$ Just as in the one-dimensional case these values retain information about initial conditions and can be found with enough data. The alternative is to use directly the initial conditions. 

\begin{figure}[H]
    \centering
    \includegraphics[width=\textwidth]{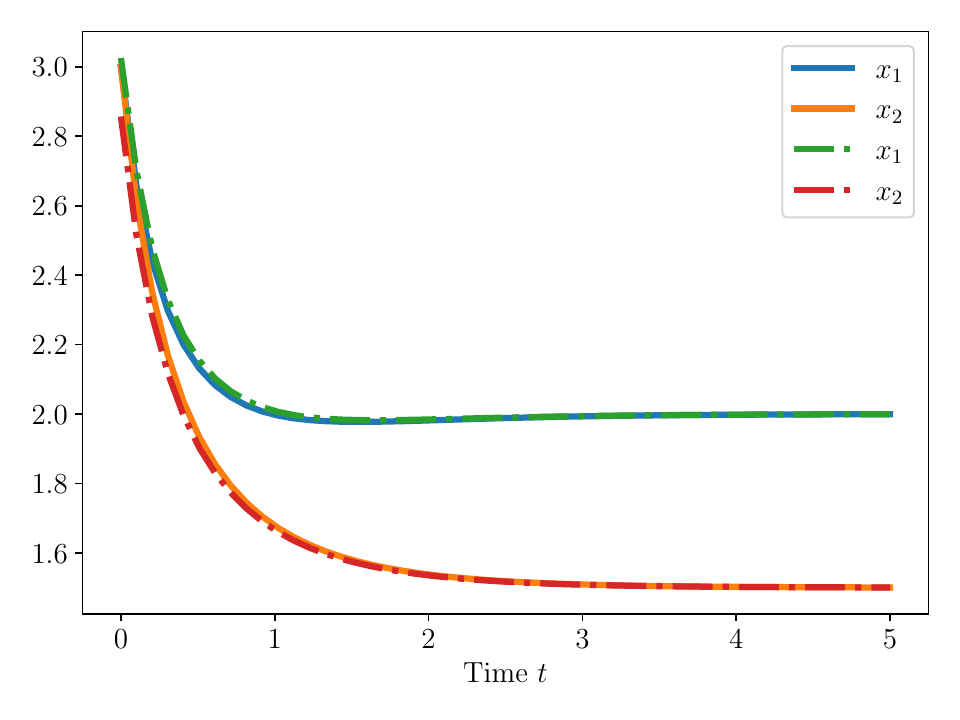}
    \caption{SPS (dashed line) and numerical (solid line) solutions, $\bm x(0) = (1,1)$. Convergence is guaranteed for $t > t_0 = 1.3$.}
    \label{fig:xy-lim}
\end{figure}

With these values, we plot our analytic SPS answer for $x_1(t)$ and $x_2(t)$ versus the image using a numerical solver (see Fig.~\ref{fig:xy-lim}).

\section{$M \times M$ systems}\label{sec:nbyn}
Now consider the $M \times M$ system:
\[ \dot{\bm x}  = \diag({\bm x}) (\bm b + \bm A \bm x) .\]
Note that existence and uniqueness of the solution is guaranteed, for any $\bm x(0) \in \mathbb{R}^M$, since the RHS is continuous and all its partial derivatives are continuous in $\mathbb{R}^M$~\cite{strogatz2018nonlinear}. In what follows we assume that the matrix $\bm A$ is invertible. Equilibrium is given by $\bm c = -\bm A^{-1} \bm b$, and again, the eigenvalues of the linearized system around equilibrium are those of the matrix $\diag(\bm c)\bm A$. 

 Let $\bm \lambda$ be a vector whose entries are the eigenvalues of $\diag(\bm c) \bm A.$ The order of the eigenvalues does not matter, but we can assume they are in ascending order of magnitude. 
We also assume that the $M-1$ dimensional hyperplane defined by $\bm z \cdot \bm \lambda = 0$ has no solutions with $\bm z$ having integer entries. (In other words, the hyperplane does not intersect the integer lattice.) With this assumption we let
the solution for each $x_i$ be of the form
\[x_i(t) = \sum_{\bm n\geq 0} \al^{\bm n}_i e^{(\bm n \cdot \bm \lambda ) t}.\]

Setting $\bm n = 0$ yields the system $\bm A \bm{\al^0} = -\bm b$, or, $\bm{\al^0} = \bm A^{-1} \bm b$, i.e., we recover that the constant terms are equal to equilibrium.

Now let us consider the case when 
$\bm n = \bm{e}_{i}$ where $\bm{e}_{i}$ is a standard basis vector with $1$ in the $i$th entry and $0$s elsewhere. 
Then $\bm n \cdot \bm \lambda$ is the eigenvalue $\lambda_{i}.$ By again equating coefficients, for these $\bm n$ our system reduces to

\[ \left ( \lambda_{i}\bm I - \diag({\bm c}) \bm A \right) \bm{\al}^{\bm n} = 0.\]

The eigenvalues must be distinct, and so it follows that the matrix
\[ \lambda_{i}\bm I - \diag({\bm c}) \bm A\] has rank $M-1.$ Thus the vector $\bm {\al}^{\bm n}$ is in the kernel of this matrix; the kernel is one dimensional. One of the coefficients must be nonzero and arbitrary and all others are determined once this value is chosen. Note also that if our initial conditions do not take on equilibrium values and in particular $x_1$, then in fact $\al_{1}^{\bm e_i}$ must not be zero. 

Thus we have shown that for the standard basis vectors, we have $M$ arbitrary values (one for each $i$). Fixing these values corresponds to choosing initial conditions for the systems.

For other values of $\bm n$, we find the following system holds:
\[\left(-\diag({\bm c}) \bm A + (\bm n \cdot \bm \lambda) \bm I \right)\bm {\al^n}= \bm {s^n},\]
where $(\bm s^n)_i = \sum_j a_{i,j} S_{i,j}^{\bm n}$ and 
\[S_{i,j}^{\bm n}= \sum_{(i_1,\dots,i_M) \neq 0, \bf{n}} \al_{i}^{i_1, i_2, \dots,i_M}  \al_{j}^{n_1 - i_1, n_2 - i_2, \dots, n_M - i_M}.\]
Let's rewrite the system as
\[(\bm n \cdot \bm \lambda) \left ( I - \frac{\diag{(\bm c)} \bm A}{\bm n \cdot \bm \lambda } \right) \bm{\al^n} = \bm{s^n}.\]

Define as before
\[ \bm S^{\bm n} = (S_{i,j}^{\bm n}) \,\,\,\,\, 1 \leq i,j \leq M .\]
The matrix $\bm S^{\bm n}$ is symmetric and it follows as in section 2 that 
$\bm s^{\bm n}$ equals the diagonal entries of $\bm A \bm {S^n}.$
Note that the matrix $\bm{n}\cdot \bm{\lambda} I - \diag({\bm c}) \bm A$ is invertible except in the case when $\bm n$ is one of the standard basis vectors. 
Thus with this set up---except for adjusting for some constants---everything done in the $2 \times 2$ case goes through.  Summarizing we have
\begin{theorem}
Suppose that the matrix $\bm A$ is invertible. Let $\bm b$ be some given vector and $\bm c = -\bm A^{-1} \bm b.$ Let $\bm \lambda$ be a vector constructed from the eigenvalues of $\diag(\bm c)\bm A.$ We assume the eigenvalues are real and negative and that there are no vector integer solutions of $\bm z \cdot \bm \lambda = 0$
Then there exists a $t_{0}$ such that
the system 
 \[ \dot{\bm x}  = \diag{(\bm x)} (\bm b + A \bm x) .\]
 has a solution of the from
 \[ \bm x = \sum_{\bm n} \bm{\al}^{n} e^{\bm n \cdot \bm{\lambda} t}\]
 that converges for  $ t \geq t_{0} >  0.$

\end{theorem}
\proof
The definitions, lemmas and conclusions from the $2 \times 2$ case extend to the higher dimensional case without almost any change except in two places. The careful reader will notice that in Lemma \ref{S.est} the conclusion 
 \[|S_{1,1}^{n_{1},n_{2}}| \leq \frac{4 K^{n_{1}+n_{2}}((\log (n_{1}+1)+1)((\log(n_{2}+1)+1)}{(n_{1}+1)(n_{2}+1)}\]
is replaced by
 \[|S_{1,1}^{\bm n}| \leq \frac{2^{M}K^{N}((\log (n_{1}+1)+1)((\log (n_{2}+1)+1) \cdots (\log (n_{M}+1)+1)}{(n_{1}+1)(n_{2}+1)\cdots (n_{M}+1)},\] where
 $N = n_{1} + n_{2} +\cdots + n_{M}.$

The other adjustment is that the Hilbert-Schmidt  norm of the matrix $\bm S^{\bm n}$ is then bounded by 
\[ \frac{2^{M}M K^{N}((\log (N_{2}+2))+1)^{M}}{(n_{1}+1)(n_{2}+1)\cdots (n_{M}+1)}.\] 

\endproof

\section{Examples}
Here are some examples that illustrate the nature and convergence of the SPS solutions.

\subsection{$2 \times 2$ system}
Let  $$ \bm A = \twotwo{-2}{-1}{-1}{-1}, \,\,\,\,\,\bm b = \left(\begin{array}{c}4\\3\end{array}\right).$$
Then $\bm c = \left(\begin{array}{c}1\\2\end{array}\right)$. The eigenvalues are given by 
$$\lambda_1 = -2 +\sqrt{2} \,\,\,\mbox{and}\,\,\, \lambda_2 = -2 -\sqrt{2}.$$
We use the solution from Sec.~\ref{sec:2by2}, and let $n_1 \leq 3$ and $n_2 \leq 3$ so that $N=6$. The coefficients $\bm \al$ in terms of $p = \al_1^{1,0}$ and $q = \al_1^{0,1}$ are given in the matrix below (to the hundredths):
\[\begin{pmatrix}
     1 & p & -0.08 p^2 & -0.64 p^3\\
 q & 2 pq & 1.28 p^2 q & -0.56 p^3 q\\
 0.93 q^2 & 2.82 p q^2&  3.51 p^2 q^2 
&  1.33 p^3 q^2\\
 0.85 q^3 & 3.46 p q^3 & 6.27 p^2 q^3 &
  5.58 p^3 q^3\end{pmatrix},
\]
where $\al_1^{i,j}$ is the $ij$-th entry, indexed from 0.

At this point, $\al_1^{1,0}$ and $\al_1^{0,1}$ are still undetermined: we must incorporate the initial condition. 
Supposing the series converge at $t = 0$, the sum of all coefficients (for one variable) must equal its initial value:
\begin{align*}
    \sum_{n_1, n_2} \al_1^{n_1,n_2} &= x_1(0)\\
\sum_{n_1, n_2} \al_2^{n_1,n_2} &= x_2(0).\end{align*}
(Note that if we do not anticipate convergence at $t = 0,$ we can use initial values at some other fixed point.)

\begin{figure}[H]
    \centering
    \begin{subfigure}[b]{0.45\textwidth}
    \centering
    \includegraphics[width=\textwidth]{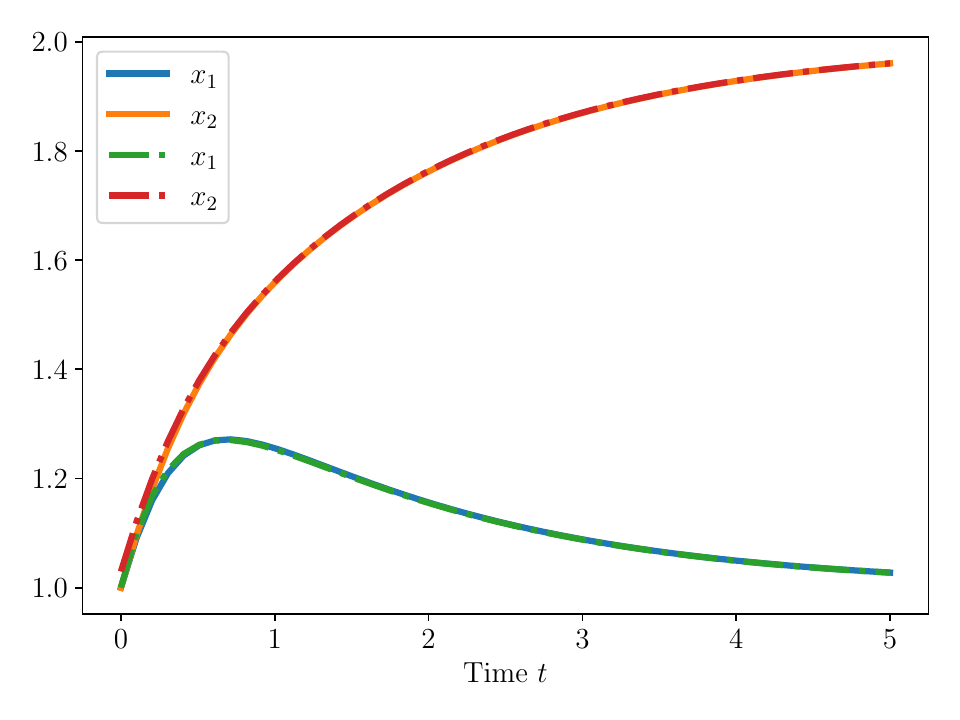}
    \caption{$\bm x(0) = (1,1)$}
    \label{fig:x1x2-11}
    \end{subfigure}
        \begin{subfigure}[b]{0.45\textwidth}
        \centering
        \includegraphics[width=\textwidth]{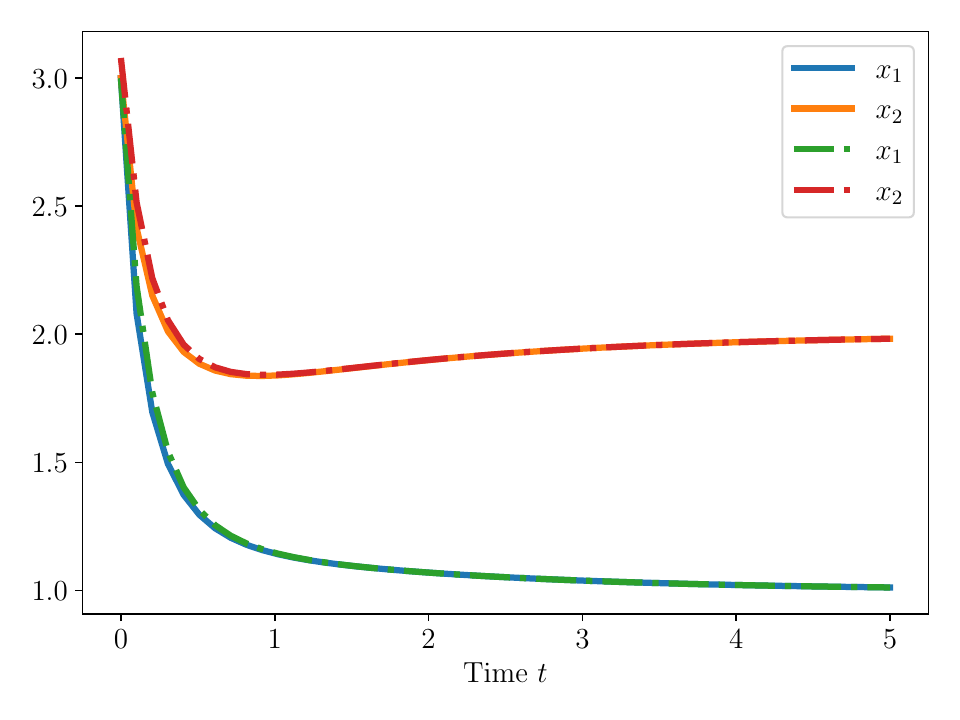}
    \caption{$\bm x(0)  = (3,3)$}
    \label{fig:x1x2-33}
    \end{subfigure}\\
    \begin{subfigure}[b]{0.45\textwidth}
        \centering
        \includegraphics[width=\textwidth]{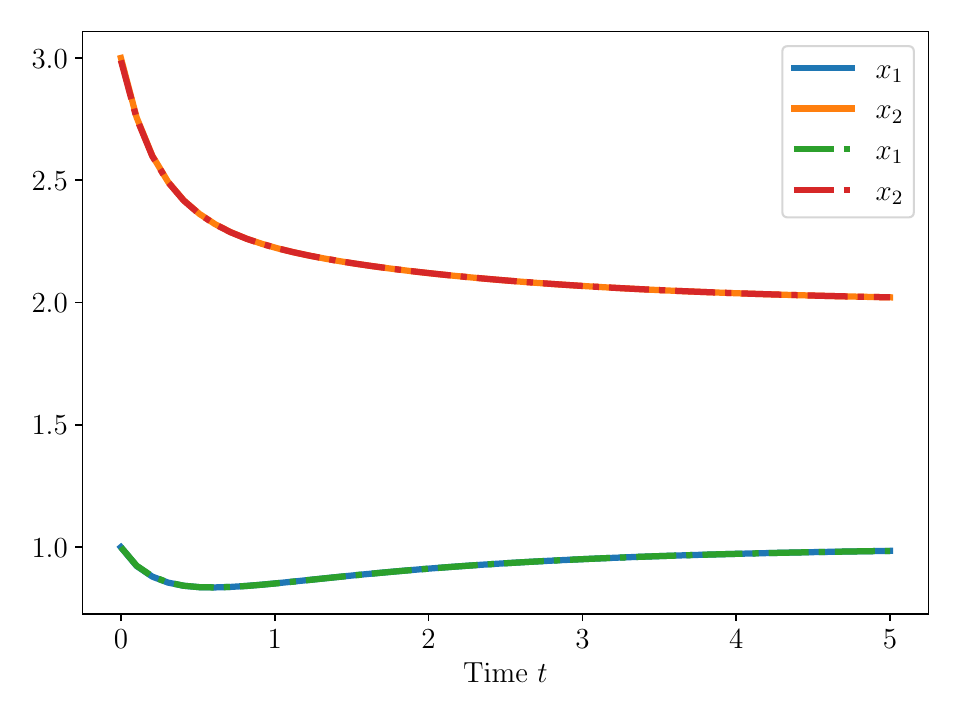}
    \caption{$\bm x(0)  = (1,3)$}
    \label{fig:x1x2-13}
    \end{subfigure}    
        \begin{subfigure}[b]{0.45\textwidth}
        \centering
        \includegraphics[width=\textwidth]{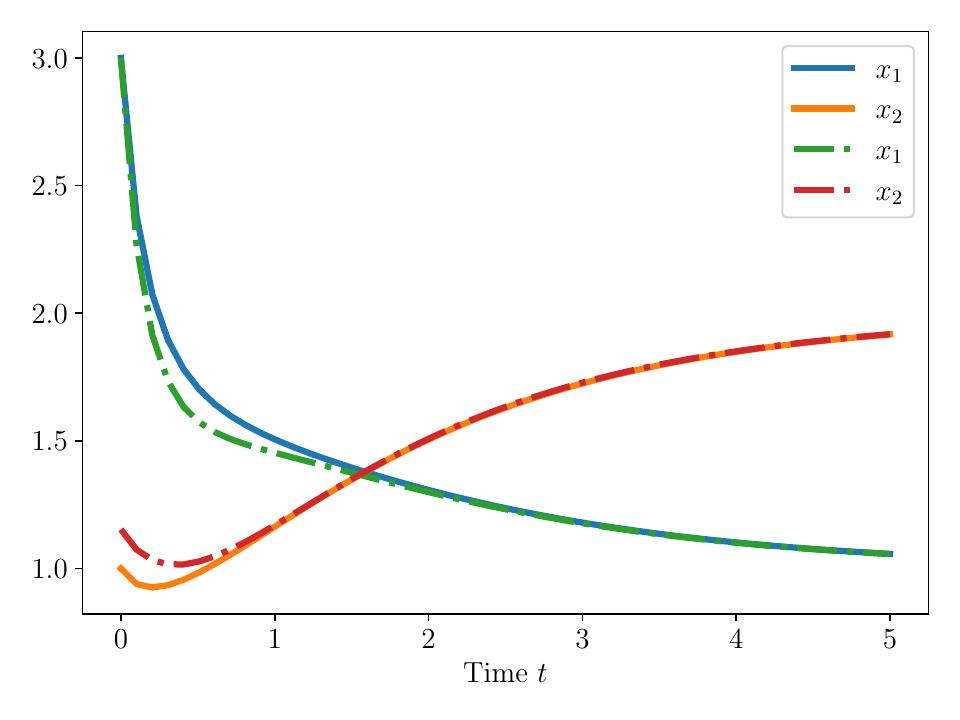}
    \caption{$\bm x(0)  = (3,1)$}
    \label{fig:x1x2-31}
    \end{subfigure} 
     \caption{Analytic SPS (dotted line) and numerical (solid line) solution trajectories for the $2 \times 2$ example, for four different initial conditions. For (A-C), analytic solutions closely match the numerical solutions for $t>0$; for (D), the analytic solution converges to the true solution around $t=2$.}\label{fig:2x2}    \end{figure}

Fig.~\ref{fig:2x2} shows the solution trajectories for four different initial conditions.

\subsection{$3 \times 3$ system}\label{ssec:ex3by3}
\begin{figure}[H]
    \centering
    \begin{subfigure}[b]{0.45\textwidth}
    \centering
    \includegraphics[width=\textwidth]{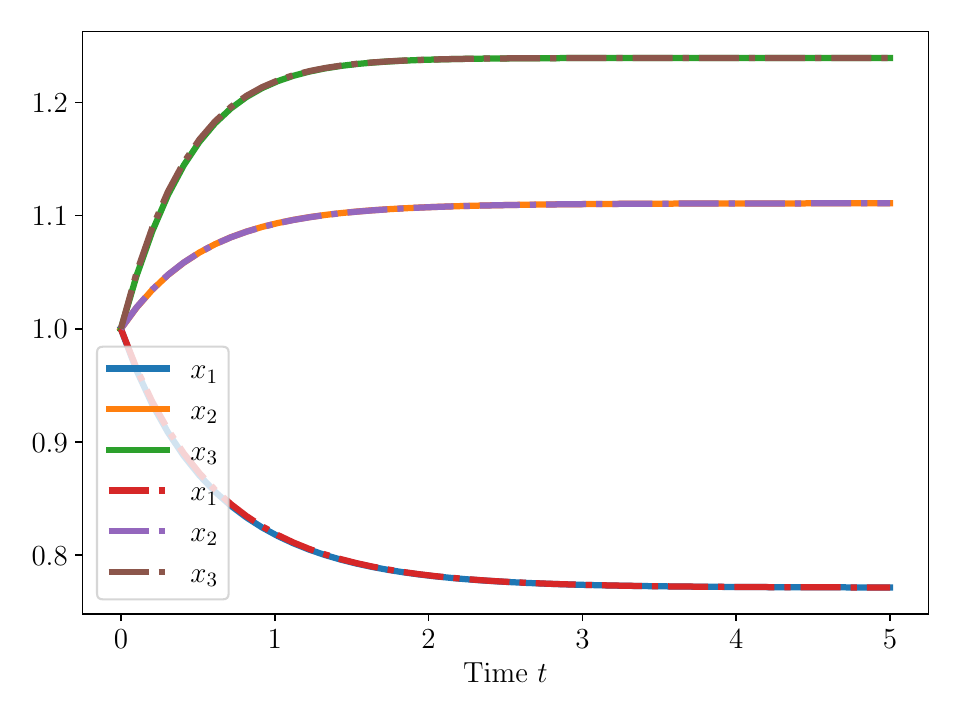}
    \caption{$\bm x(0)  = (1,1,1)$}
    \label{fig:xyz2-111}
    \end{subfigure}
        \begin{subfigure}[b]{0.45\textwidth}
        \centering
        \includegraphics[width=\textwidth]{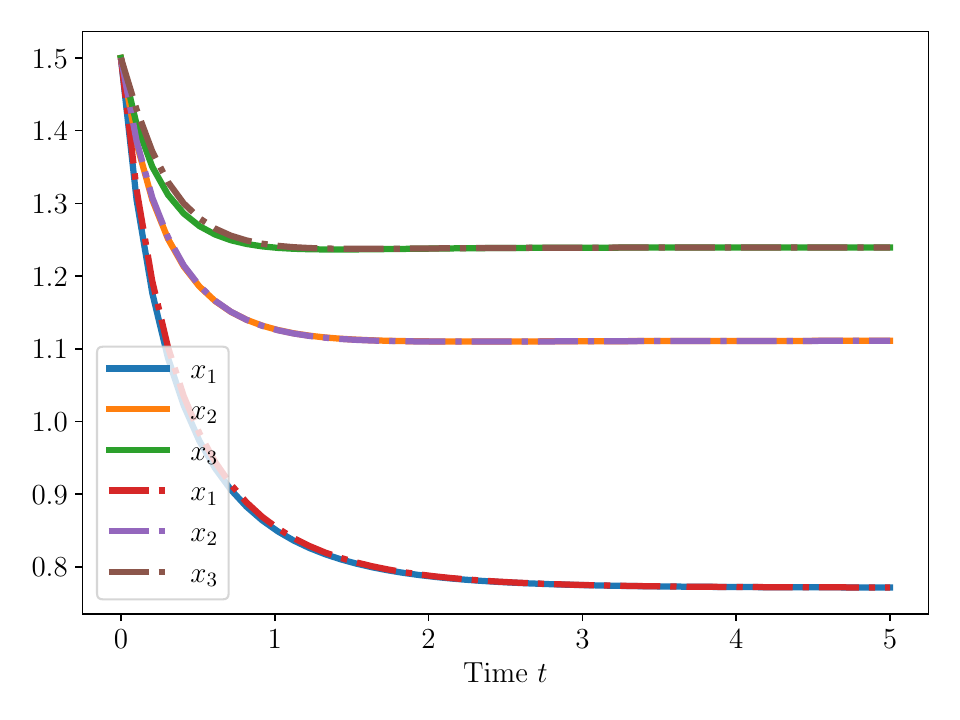}
    \caption{$\bm x(0)  = (3/2,3/2,3/2)$}
    \label{fig:xyz2-1p5}
    \end{subfigure}\\
    \begin{subfigure}[b]{0.45\textwidth}
        \centering
        \includegraphics[width=\textwidth]{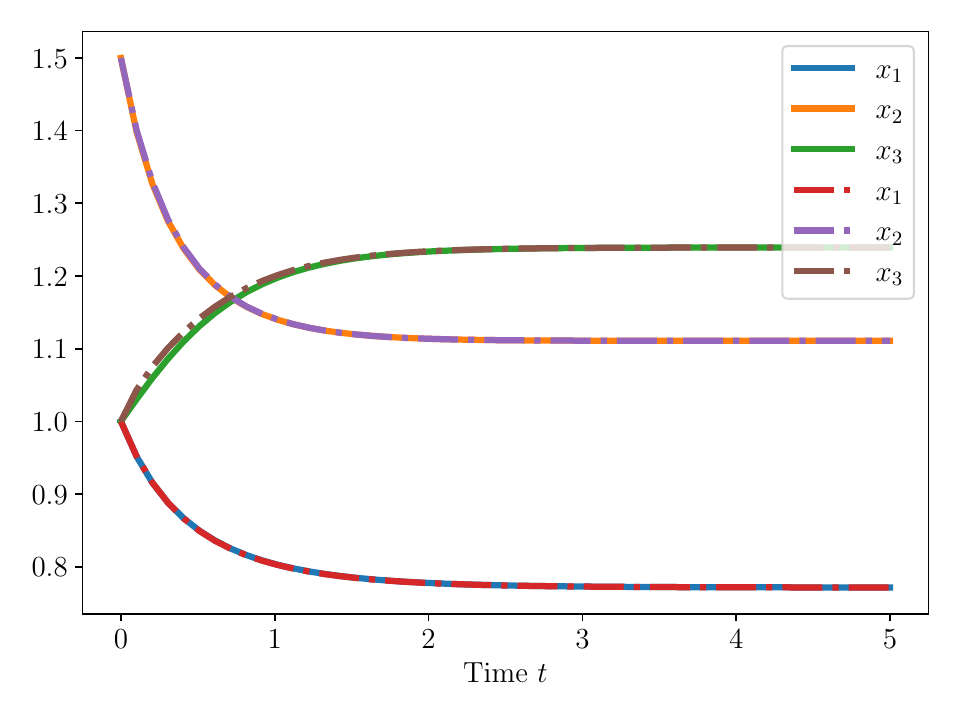}
    \caption{$\bm x(0)  = (1,3/2,1)$}
    \label{fig:xyz2-ic3}
    \end{subfigure}    
        \begin{subfigure}[b]{0.45\textwidth}
        \centering
        \includegraphics[width=\textwidth]{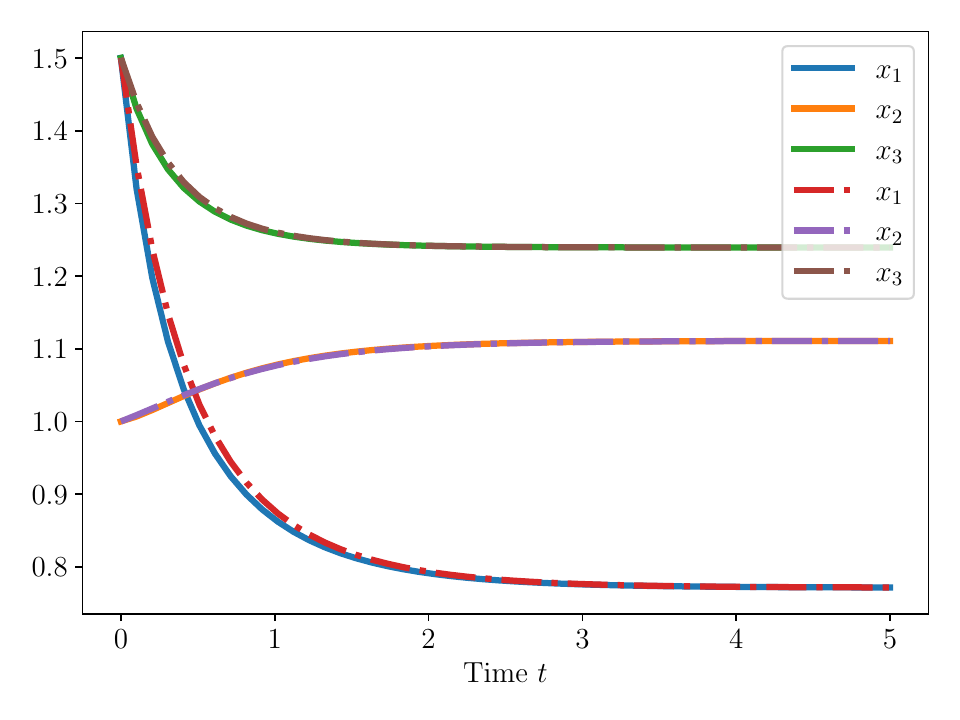}
    \caption{$\bm x(0)  = (3/2,1,3/2)$}
    \label{fig:xyz2-ic4}
    \end{subfigure} 
     \caption{Analytic SPS (dotted line) and numerical (solid line) solution trajectories for the $3 \times 3$ example, for four different initial conditions. Analytic solutions closely match the numerical solutions for $t>0$.}\label{fig:3x3}    \end{figure}

Let's look at another example, now a $3 \times 3$ system. Let the interaction matrix and linear growth rate vector be
\[ \bm A = \begin{pmatrix}
    -2 & -0.3 & -0.1 \\ -0.2 & -2 & -0.1 \\ -0.1 & -0.4 & -2
\end{pmatrix}, 
\quad \bm b = \begin{pmatrix}
    2 \\ 2.5 \\ 3
\end{pmatrix}.\]
Again we let $n_1, n_2, n_3 \leq 3$ so that $N = 9$. Fig.~\ref{fig:3x3} shows numerical and analytic solution trajectories for the three variables, for four different initial conditions.

\section{Remarks about Koopman expansions}

In~\cite{mezic2019} (see displayed equation 69), a general framework is given for the systems we have considered and a Koopman expansion of the form 
\[ \bm x (t) = \sum_{\bm n} \bm{\al}^{\bm n} e^{\bm n \cdot \bm{\lambda} t}\] 
is produced. Can we see that they are the same?

Let us consider the $2 \times 2$ system. Since the SPS coefficients are all generated by products of $\al^{10}_{1}$ and $\al^{01}_{1}$ and these values can be found by either a limiting process or an algebraic method, if the series that appears in the above reference is known to converge then the coefficients must be the same. 

But we can say more. Let us make the assumption that $\al^{10}_{1}$ and $\al^{01}_{1}$ are eigenfunctions of the Lie derivative operator for the Koopman flow with eigenvalues $\lambda_{1}$ and $\lambda_{2}.$ 
Recall our system is of the form:
\begin{align*} \dot{x_1} &= f(x_1,x_2) =b_1x_1 + a_{11} x_1^{2} + a_{12} x_1 x_2 \\
\dot{x_2} &= g(x_1,x_2) = b_2 x_2  + a_{21} x_1x_2 + a_{22} x_2^2.\end{align*} 
For $\al^{10}_{1}$ and $\al^{01}_{1}$ to be eigenfunctions it follows that
  $$\lambda_{1}\al^{10}_{1} = f(x_{1}, x_{2})\frac{\partial\, \al^{10}_{1}}{\partial x_1} + g(x_{1}, x_{2})\frac{\partial\, \al^{10}_{1}}{\partial x_2}$$ and
   $$\lambda_{2}\al^{01}_{1} = f(x_{1}, x_{2})\frac{\partial\, \al^{01}_{1}}{\partial x_1}+ g(x_{1}, x_{2})\frac{\partial\, \al^{01}_{1}}{\partial x_2}.$$
  Then it is the case that $\al^{\bm{n}}_{i}$ does align as well with eigenvalue $\bm n \cdot \bm{\lambda}.$ This is easy to see by recalling from Theorem 2.6 that all
$\al_i^{n_1,n_2}, i = 1,2$ have the form 
\[ c_{\bm n} (\al_1^{1,0})^{n_{1}} (\al_1^{0,1})^{n_{2}} \] where the constant 
$c_{\bm n}$ only depends on all the other data in the system. 

We now show that the functions
$\al^{10}_1$ and $\al^{01}_1$ satisfy the Koopman eigenfunction conditions.

Recall that our solution is of the form
\[ \bm x(t) =  \sum_{\bm n}c_{\bm n} (\al_1^{1,0})^{n_{1}} (\al_1^{0,1})^{n_{2}}e^{\bm \lambda \cdot \bm n \,t},\] and thus
\[ \dot{\bm x}(t) =  \sum_{\bm n}c_{\bm n} (\al_1^{1,0})^{n_{1}} (\al_1^{0,1})^{n_{2}} (\bm \lambda \cdot \bm n )\,e^{\bm \lambda \cdot \bm n \,t}.\]
Assuming we have convergence at $t = 0,$ this means that 
\begin{equation}\label{x1}  x_1(0)  =  \sum_{\bm n}c_{1,\bm n} (\al_1^{1,0})^{n_{1}} (\al_1^{0,1})^{n_{2}}\end{equation}
and 
\begin{equation}\label{x2} x_2(0)  =  \sum_{\bm n}c_{2,\bm n} (\al_1^{1,0})^{n_{1}} (\al_1^{0,1})^{n_{2}}.\end{equation}
Now think of $x_1(0) = x_1$ and $x_2(0) = x_2$ as variables. If we take the partial derivative of 
$$\sum_{\bm n}c_{1,\bm n} (\al_1^{1,0})^{n_{1}} (\al_1^{0,1})^{n_{2}}$$ with respect to $x_1$ it is a straightforward computation to see that it is
$$\sum_{\bm n}c_{1,\bm n} (\al_1^{1,0})^{n_{1}} (\al_1^{0,1})^{n_{2}}\left(n_1 \frac{\frac{\partial \al_1^{10}}{\partial x_1}}{\al_1^{10}} +  n_2 \frac{\frac{\partial \al_1^{01}}{\partial x_1}}{\al_1^{01}} \right).$$
The equivalent statement can be made for the above with $c_{1, \bm n}$ replace with $c_{2, \bm n}$
and also for the $\partial$s with respect to $x_1$ replaced by $x_2.$ 

In what follows we define
$k_1 = \ln \al_1^{10}\,\,\, k_2 = \ln \al_1^{01}$ and note that 
$$\frac{\partial k_1}{\partial x_1} = (k_1)_{x_1} = \frac{\frac{\partial \al_1^{10}}{\partial x_1}}{\al_1^{10}},\,\,\,\, \,\,\,\,\frac{\partial k_1}{\partial x_2} = (k_1)_{x_2} = \frac{\frac{\partial \al_1^{10}}{\partial x_2}}{\al_1^{10}}, $$
and
$$\frac{\partial k_2}{\partial x_1} = (k_2)_{x_1} = \frac{\frac{\partial \al_1^{01}}{\partial x_1}}{\al_1^{01}},\,\,\,\, \,\,\,\,\frac{\partial k_2}{\partial x_2} = (k_2)_{x_2} = \frac{\frac{\partial \al_1^{01}}{\partial x_2}}{\al_1^{01}}. $$
Now define a matrix from these partials:
$$K = \twotwo{(k_1)_{x_1}}{(k_2)_{x_1}}{(k_1)_{x_2}}{(k_2)_{x_2}}$$ and two vectors 
$$\bm v_1 = \left(\begin{array}{c} \sum_{\bm n}c_{1,\bm n} (\al_1^{1,0})^{n_{1}}(\al_1^{0,1})^{n_{2}}\,n_1 \\
 \sum_{\bm n}c_{1,\bm n} (\al_1^{1,0})^{n_{2}}(\al_1^{0,1})^{n_{2}}\,n_2
\end{array}\right)$$
$$\bm v_2 = \left(\begin{array}{c}  \sum_{\bm n}c_{2,\bm n} (\al_1^{1,0})^{n_{1}}(\al_1^{0,1})^{n_{2}}\,n_1 \\
 \sum_{\bm n}c_{2,\bm n} (\al_1^{1,0})^{n_{2}}(\al_1^{0,1})^{n_{2}}\,n_2
\end{array}\right).$$

With these definitions and taking the partials of both sides of the equations in (\ref{x1}) and (\ref{x2}) we have that
\[ K \,\bm v_1 = \left(\begin{array}{c}1 \\0\end{array}\right),\,\,\, K\,\bm v_2 = \left(\begin{array}{c}0 \\1\end{array}\right). \]
We also have that 
\[ \bm v_1 \cdot \left(\begin{array}{c}\lambda_1 \\\lambda_2\end{array}\right) = f\]
where $f$ corresponds to the first row (on the right) of our system, and
\[ \bm v_2 \cdot \left(\begin{array}{c}\lambda_1 \\\lambda_2\end{array}\right) = g\] with $g$ being the second row. 
Thus
\[ \langle K^{-1} \left(\begin{array}{c}1 \\0\end{array}\right) , \bm \lambda \rangle  = f,\,\,
\langle K^{-1} \left(\begin{array}{c}0 \\1\end{array}\right) , \bm \lambda \rangle  =  g \]or
\[ \langle  \left(\begin{array}{c}1 \\0\end{array}\right) , (K^{-1})^t\bm \lambda \rangle  = f,\,\,
\langle  \left(\begin{array}{c}0 \\1\end{array}\right) , (K^{-1})^{t}\bm \lambda \rangle  =  g .\]
This means that 
\[ (K^{-1})^{t}\bm \lambda  = \left(\begin{array}{c}f \\g\end{array}\right),\]and the result immediately follows from this.

This is analogous to the one-dimensional case where we chose $\frac{x_{0} -k}{x_{0}}$ to be our parametrization. The higher-dimensional cases can be analyzed in a similar fashion.
To summarize we have shown that our solution is a convergent (for large enough $t$) Koopman expansion with explicitly constructed coefficients. 

\section{Reduced models}
In previous work~\cite{morrison2020data, morrison2020embedded}, we have seen that reduced models, that only track a subset of $L<M$ variables, can still faithfully represent the dynamics of the complete $M\times M$ system.
We lacked a theoretical justification for why this should be the case, but now the analysis here can help explain this phenomenon.

In practice, this type of reduction is very common. In hydrocarbon combustion, a state-of-the-art model may include fifty (or more) chemical species, yielding a set of fifty ODEs~\cite{grimech}. But these models are too expensive for many industrial and academic pursuits; instead reduced mechanisms include only a handful of chemical species, such as just the reactants and products. As another example, consider SEIR-type models of epidemiological outbreaks (which track susceptible, exposed, infected, and recovered populations, or other groups such as hospitalized or quarantined). But of course many other variables (epidemiological species) may contribute to the outbreak dynamics. Typical models of a Zika virus outbreak track sub-populations of humans and mosquitoes~\cite{dantas2018calibration}, but omit livestock, non-human primates, and so on. Modelers may not have sufficient data to include these variables, or simply not know which are the most relevant.

Despite such omission of potentially many state variables, still reduced models may capture observed dynamics and even perform well when extrapolated to other scenarios (such as in the future, or for other environmental conditions)~\cite{morrison2020embedded}. To examine this type of reduced model over $L<M$ state variables, let us first take the \emph{partial} model, using $A_{L\times L}$, the upper left $L\times L$ submatrix of $A$ and $b_{1:L}$, the first $L$ components of $b$.

Following~\cite{morrison2020data, morrison2020embedded} and as a correction to the partial model, two terms are added for each RHS:
\[\dot{x_i} = b_i x_i + (a_{i,1} x_1 + \dots + a_{i,L} x_L)x_i + \delta_{i}x_i + \gamma_{i}\dot x_i .\]
For each $i$, $\delta_{i}$ and $\gamma_{i}$ are calibrated using data from the true (complete) system. For details of the calibration process, see~\cite{morrison2020embedded}.

It is easy to see that $\bm \delta$ corrects for equilibrium and $\bm \gamma$ adjusts the overall transient behavior. Let $(\gamma^*)_i = \frac{1}{1-\gamma_{i}}$.
We can rewrite the system with
\[\hat{\bm A} = \diag{(\bm\gamma^*)} \bm A_{L \times L},
\quad \quad \hat{\bm b} = \diag(\bm \gamma^*) (\bm b_{1:L} + \bm\delta)  .\]
Then the equilibrium of this system is given by
\[\hat{\bm c} = -\bm A^{-1} (\bm b_{1:L} + \bm \delta)\]
and the eigenvalues of the new linearized system are the eigenvalues of
\[\diag{(\bm \gamma^*)} \diag(\hat{\bm c}) \bm A_{L\times L}.\]

So we conjecture that a good reduced model sets $\bm \delta$ such that $\hat{\bm c} = \bm c_{1:L}$ (i.e., set $\bm \delta = -(\bm A \bm c + \bm b)_{1:L}$)
 and sets $\bm \gamma$ so that the $L$ eigenvalues (in decreasing order) match the first $L$ eigenvalues of the $\diag{(\bm c)} \bm A$. Since we have $L$ unknowns (degrees of freedom) in the vector $\bm \gamma$ and we have a target of the first $L$ eigenvalues this should be possible algebraically, at least in some generic sense. 

\begin{figure}[H]
    \centering
    \begin{subfigure}[b]{0.45\textwidth}
    \centering
    \includegraphics[width=\textwidth]{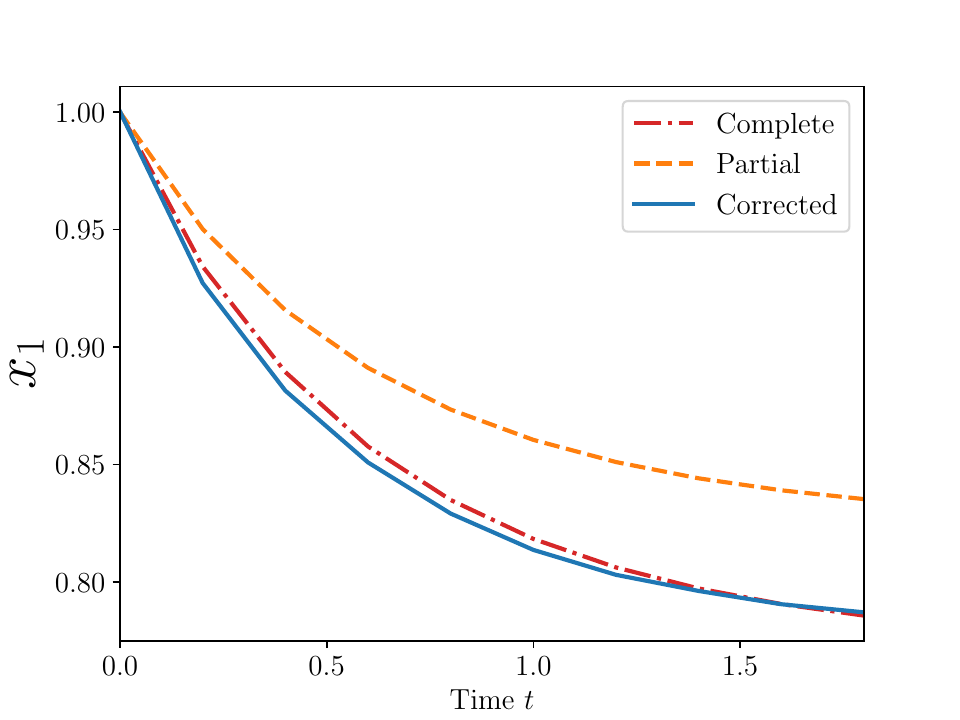}
    \caption{$\bm x(0) = (1,1,1) $}
    \label{fig:x1-red}
    \end{subfigure}
        \begin{subfigure}[b]{0.45\textwidth}
        \centering
        \includegraphics[width=\textwidth]{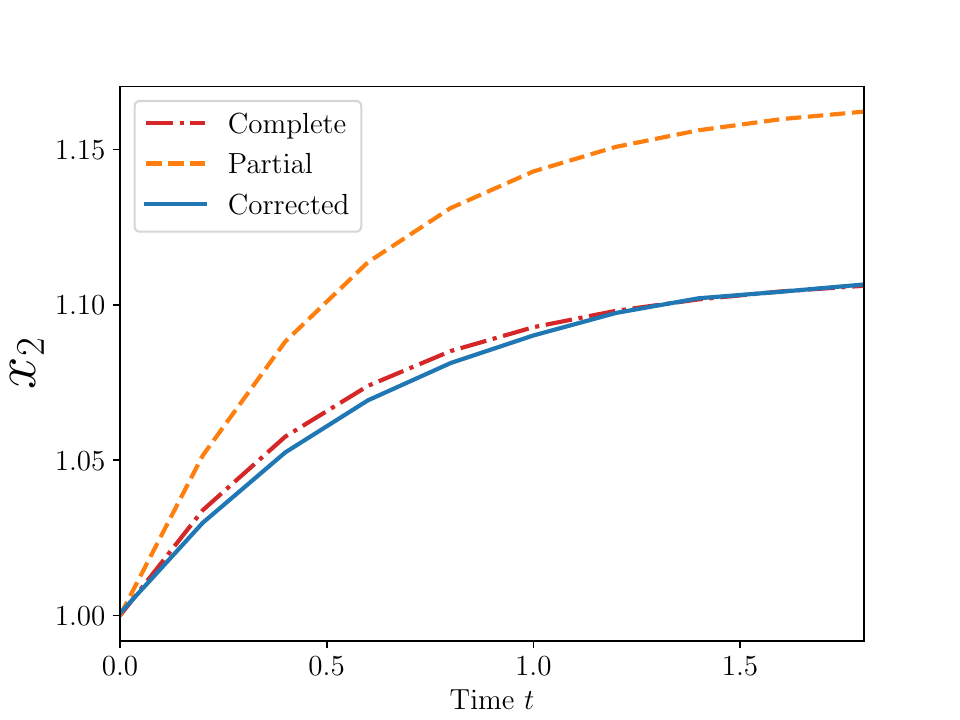}
    \caption{$ \bm x(0) = (1,1,1)$}
    \label{fig:x2-red}
    \end{subfigure}
        \begin{subfigure}[b]{0.45\textwidth}
    \centering
    \includegraphics[width=\textwidth]{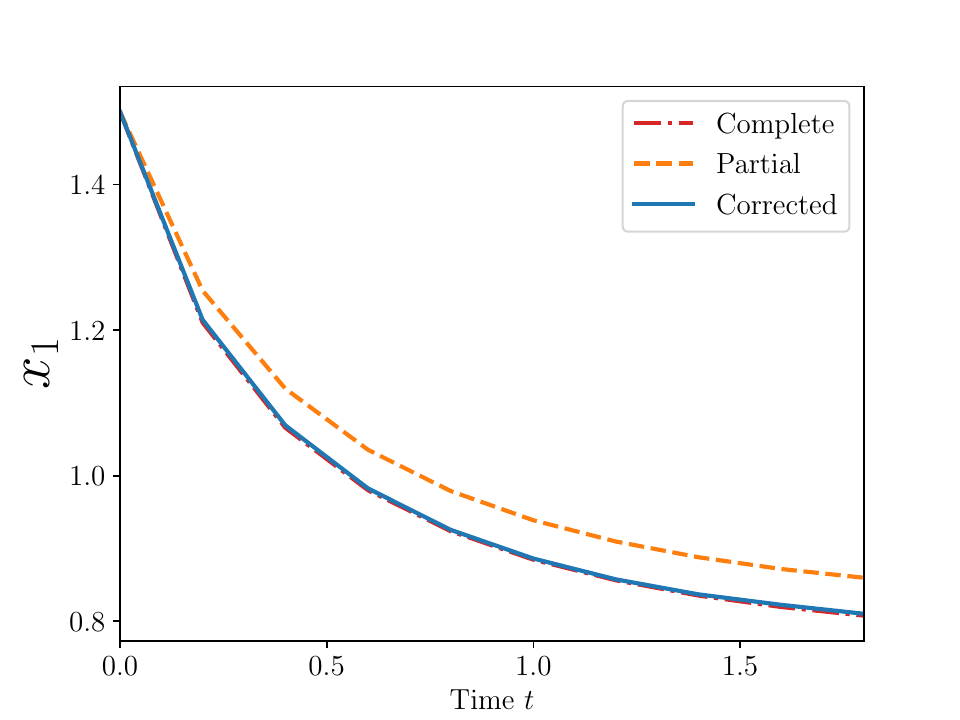}
    \caption{$ \bm x(0) = (3/2,1/2,3/2)$}
    \label{fig:x1-red2}
    \end{subfigure}
        \begin{subfigure}[b]{0.45\textwidth}
        \centering
        \includegraphics[width=\textwidth]{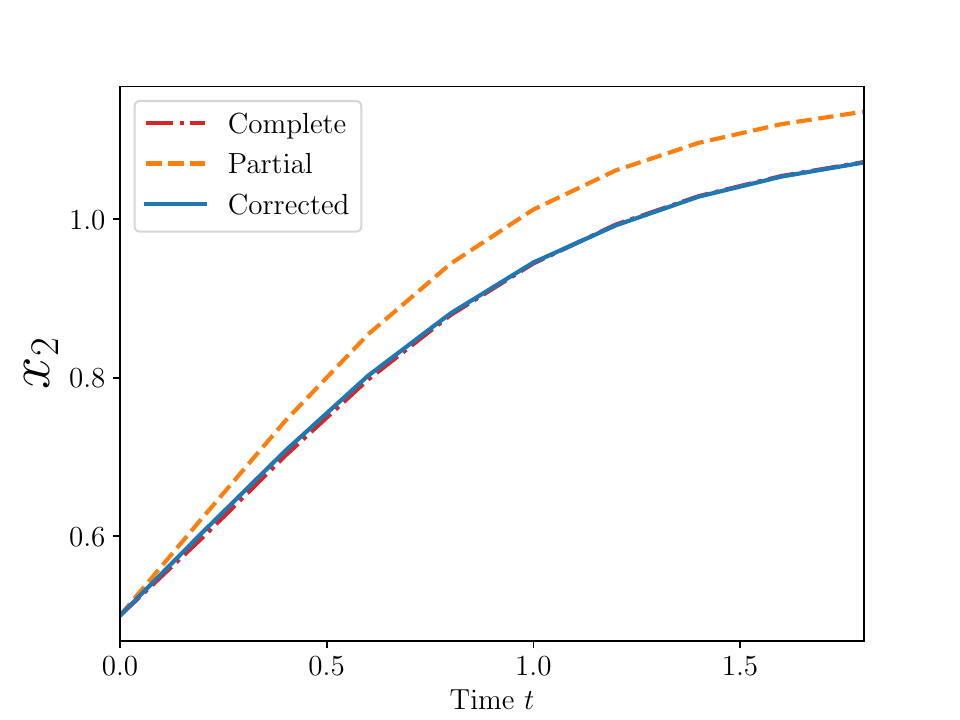}
    \caption{$\bm x(0) = (3/2,1/2,3/2)$}
    \label{fig:x2-red2}
    \end{subfigure}
    \caption{Numerical solutions for complete, partial, and corrected models. The complete model includes $3$ variables, but the corrected model only includes $2$ of those. The partial model retains the relevant sub-matrix and sub-vector from the complete model.}\label{fig:xy-red}
\end{figure}

As an example, let's consider the system in section~\ref{ssec:ex3by3}, and reduce the system to 2 variables:
\[ \bm A_{L \times L} = \begin{pmatrix}
    -2 & -3/10 \\ -2/10 & -2
\end{pmatrix}, 
\quad \bm b_{1:L} = \begin{pmatrix}
    2 \\ 5/2
\end{pmatrix}.\]
In this example, we find that for $x_1$, $\delta_{1} \approx -0.14$ and $\gamma_{1} \approx -0.01$, and for $x_2$, $\delta_{2} \approx -0.11$ and $\gamma_{2} \approx -0.07$.
In fact, this is what happens for this example. With the values of $\bm \delta$ and $\bm \gamma$ given above, we find the following, to the hundredths:
\begin{itemize}
    \item $c_1 = 0.77$, $\hat c_1 = 0.76$
    \item $c_2 = 1.11$, $\hat c_2 = 1.12$
    \item $\lambda_1 = -1.48$, $\hat \lambda_1 = -1.44$
    \item $\lambda_2 = -2.11$, $\hat \lambda_2 = -2.16$.
\end{itemize}
Finally, note that the four values above of equilibrium and eigenvalues ($c_1$, $c_2$, $\lambda_1$, $\lambda_2$) do not depend on the initial conditions, so a corrected reduced model may be valid over a range of initial conditions, and not just for a single trajectory. This broad domain of validity is critical if the reduced model is meant to be used in other contexts beyond the calibration scenario, which is typically what one asks of a model (Fig.~\ref{fig:xy-red}).

\section{Discussion}
Analytic solutions for nonlinear sets of differential equations are few and far between. We find analytic solutions to a class of nonlinear coupled ODES with quadratic RHS via power series expansions, and furthermore, we:
\begin{enumerate}
    \item Show, for the logistic equation, if \emph{any} power series solution exists, it must be of the exponential form given here;
    \item Derive an explicit construction for the coefficients and how they decay; 
    \item Prove these power series solutions converge for $t > t_0$ and give a bound;
    \item Compare the analytic solutions to numerical solutions through various experiments;
    \item Show how, considering the eigenvalues of the linearized system, reduced systems may still capture much of the system dynamics. 
\end{enumerate}

We conclude with a couple remarks about these solutions. First, this work reveals the exponential nature of the trajectories, even away from equilibrium. For a given problem and for $t > t_0$, the solutions are infinitely smooth and differentiable, and approach equilibrium exponentially fast. Besides its theoretical interest, this may inform numerical analysis issues such as convergence, order, and stability, if a numerical solver is still to be used.

Second, numerical experiments typically result in coefficients which decay quickly. Thus, a truncated series even with relatively small $N$, say $N \approx \gamma M$ where $\gamma = 2$ or $3$ can result in very good closed-form approximations to the true solution. If one simply needs a quick estimate of the solution, a truncated series comprised of only a handful terms may do the trick.
Third, nothing in this work seems to fundamentally rely on the particular form of the RHS. We propose that many sets of coupled autonomous ODES with polynomial RHS will admit a similar exponential power series solution. Of course, solving for the coefficients may become more complicated (with, e.g., higher-order interactions or complex eigenvalues) but we expect these types of solutions, at least for some semi-bounded time-domain, to exist and converge.

\section{Concluding Remarks}
The first author has admired and been inspired by the work of Albrecht B\"{o}ttcher for almost her entire career. His creative and powerful use of linear algebra operator theory to solve problems in analysis transformed the study of asymptotics of determinants of structured operators. While this paper is in a different area, we hope that it shows that using a bit of linear algebra one can produce concrete solutions to some classical problems. 
\bibliographystyle{plain}
\bibliography{refs.bib}

Estelle Basor\\
ebasor@aimath.org

Rebecca Morrison\\
rebeccam@colorado.edu
\end{document}